\documentclass[11pt,a4paper]{article}
\usepackage{tikz}
\usetikzlibrary{patterns}

\usepackage{epsfig}
\usepackage{graphics}
\usepackage{amsmath,amssymb}
\usepackage{subfigure}
\usepackage{graphicx}
\usepackage{epstopdf}
\usepackage{hyperref}
\usepackage{pdfsync}
\usepackage{color}
\usepackage{ulem}

\newtheorem{lemma}{Lemma}
\newtheorem{proposition}{Proposition}

\numberwithin{equation}{section}
\numberwithin{thm}{section}
\numberwithin{lemma}{section}
\numberwithin{proposition}{section}
\numberwithin{definition}{section}
\numberwithin{remark}{section}
\numberwithin{figure}{section}
\numberwithin{table}{section}

\title{Minimax testing in a statistical inverse problem with unknown operator}
\author{{\em Cl\'ement ~Marteau}, \\
         Univ Lyon, Universit\'e Claude Bernard Lyon 1,\\     CNRS UMR 5208, Institut Camille Jordan, \\
         43 blvd. du 11 novembre 1918, \\F-69622 Villeurbanne cedex, France\\
        {\em Email}:~\texttt{marteau@math.univ-lyon1.fr}
          \\ \\
        and
          \\ \\
        {\em Theofanis ~Sapatinas},\\
        Department of Mathematics and Statistics,\\
        University of Cyprus,\\
        P.O. Box 20537,
        CY 1678 Nicosia,
        Cyprus.\\
        {\em Email}:~\texttt{fanis@ucy.ac.cy}}

\begin{document}
\maketitle

\begin{abstract}
We study minimax testing in a statistical inverse problem when the associated operator is unknown. In particular, we consider observations from the Gaussian regression model $y_k = b_k \theta_k + \epsilon\, \xi_k$, $k \in \mathbb{N}_m:=\{1,2,\ldots,m\}$, for some $m <+\infty$, where $\theta=(\theta_k)_{k \in \mathbb{N}_m}$ is the unknown signal of interest, $b= (b_k)_{k \in \mathbb{N}_m}>0$ is an unknown collection (the singular values of the associated finite rank operator), $\xi=(\xi_k)_{k \in \mathbb{N}_m}$ is a collection of independent standard Gaussian random variables and $\epsilon>0$ is the noise level.  It is assumed that the collection $b= (b_k)_{k \in \mathbb{N}_m}$ is unknown but contained in a given dictionary $\mathcal{B}$ of finite cardinality $|\mathcal{B}|$. Using the non-asymptotic framework for minimax testing (that is, for any fixed value of the noise level $\epsilon>0$), we provide optimal separation conditions for the goodness-of-fit testing problem. 
We restrict our attention to the specific case where the dictionary contains only two members, namely $|\mathcal{B}|=2$. 
As we will demonstrate, even this simple case is quite intrigued and reveals an interesting phase transition phenomenon. The case $|\mathcal{B}| > 2$ is even more involved, requires different strategies, and it is only briefly discussed.  \\


\end{abstract}

\newpage

\section{Introduction}
We consider the following heteroscedastic statistical model (HSM)
\begin{equation}
y_k = b_k \theta_k + \epsilon\, \xi_k, \quad k\in \mathcal{I},
\label{eq:gsm}
\end{equation}
where $\mathcal{I}=\mathbb{N}_m:=\lbrace 1,2,\ldots,m \rbrace$ for some $m<+\infty$ or $\mathcal{I}=\mathbb{N}:=\lbrace 1,2,\ldots, \rbrace$, $\theta=(\theta_k)_{k \in \mathcal{I}}$ is the unknown signal of interest, $b= (b_k)_{k \in \mathcal{I}}$ is a positive and decreasing sequence, $\xi=(\xi_k)_{k \in \mathcal{I}}$ is sequence of independent standard Gaussian random variables and $\epsilon>0$ is the noise level. The observations are given by the sequence $y=(y_k)_{k \in \mathcal{I}}$ and their joint law is denoted by $\mathbb{P}_{\theta,b}$. The HSM (\ref{eq:gsm}) is alternatively called a Gaussian regression model (GRM) when $\mathcal{I}=\mathbb{N}_m$ for some $m<+\infty$,  and a Gaussian sequence model (GSM) when  $\mathcal{I}=\mathbb{N}$. The HSM (\ref{eq:gsm}) has been widely studied in the statistical literature, especially in the estimation context (see, e.g., \cite{Bissantz_2007}, \cite{Cavalier_book}, \cite{cavalier2}).  \\

The HSM (\ref{eq:gsm}) has been mainly used to describe statistical inverse problems. Indeed, for a bounded linear (injective) operator $A$ acting on a Hilbert space $H$ with inner product $\langle \cdot , \cdot \rangle_H$ with values on another Hilbert space $K$ with inner product $\langle \cdot , \cdot \rangle_K$, consider an unknown function $f \in H$ indirectly observed in the Gaussian white noise model (GWNM), i.e.,
\begin{equation}
Y(g) = \langle Af,g \rangle_K + \epsilon\, \xi(g), \quad g \in K,
\label{eq:GWNM}
\end{equation}
where $\xi(g)$ is assumed to be a centered Gaussian random variable with variance $\|g\|_K^2=\langle g,g \rangle_K$. 
If $A$ is assumed to be compact, it admits a singular value decomposition (SVD) $(b_k, \psi_k,\phi_k)_{1 \leq k \leq N}$, where either $N=m <+\infty$ (if $A$ is a finite rank operator) or $N=+\infty$.  In particular, for orthonormal eigenfunctions $(\psi_k)_{1 \leq k \leq N}$ and $(\phi_k)_{1 \leq k \leq N}$, we have
\begin{equation}
\label{eq:svd}
A\phi_k=b_k\psi_k, \quad A^\star\psi_k=b_k\phi_k, \quad k=1,2,\ldots,N,
\end{equation}
where $A^\star$ denotes the adjoint operator of $A$ and $(b_k)_{1 \leq k \leq N}>0$ are the so-called singular values of the operator $A$ . (Note that $(b^2_k)_{1 \leq k \leq N}$ and $(\phi_k)_{1 \leq k \leq N}$ are, respectively, the eigenvalues and eigenfunctions of $A^\star A$.)
Using (\ref{eq:GWNM}) and (\ref{eq:svd}), we get observations
$$
y_k := Y(\psi_k) = \langle Af, \psi_k \rangle_K + \epsilon\, \xi(\psi_k) = \langle f, A^\star \psi_k \rangle_H + \epsilon\, \xi(\psi_k) := b_k \theta_k + \epsilon\, \xi_k, \quad  k=1,2,\ldots,N,
$$
where $\theta_k:=\langle f, \phi_k \rangle$, $1 \leq k \leq N$, and the $\xi_k:=\xi(\psi_k)$, $1 \leq k \leq N$, are independent standard Gaussian random variables. Observations from the GSM are obtained when $H$ and $K$ are infinite-dimensional Hilbert spaces. In this case,  $N=\infty$ and $0<b_k \rightarrow 0$ as $k \rightarrow +\infty$.
Observations from the GRM are obtained when either $H$ or $K$ is a finite-dimensional Hilbert space so that $A$ is a finite rank operator and, hence, compact. In this case, $N={\rm rank}(A)=m$ and the $(b_k)_{1 \leq k \leq m}$ can be arranged so that $b_1 \geq b_2 \geq \ldots \geq b_m>0$.  \\

The GSM has attracted attention in the minimax testing literature. Usually, the sequence $b=(b_k)_{k \in \mathbb{N}}$ (or, equivalently, the corresponding operator $A$) is assumed to be fully known (see, e.g., \cite{ISS_2011}, \cite{ISS_2012}, \cite{IS_2003}, \cite{LLM_2011}, \cite{LLM_2012}, \cite{MS_2015}).  Some attempts have been recently proposed to relax the assumption that the sequence $b= (b_k)_{k \in \mathbb{N}}$ (or, equivalently, the corresponding operator $A$) is available, assuming that it is partially known in the sense that we have at our disposal observations on the sequence $b= (b_k)_{k \in \mathbb{N}}$ of the form $x_k = b_k+ \eta \zeta_k \quad k\in \mathbb{N}$,
where $\zeta=(\zeta_k)_{k \in \mathbb{N}}$ is a sequence of independent standard Gaussian random variables (that is also independent of the sequence $\xi=(\xi_k)_{k \in \mathbb{N}}$) and $\eta>0$ is a noise level (see, e.g., \cite{Kroll_2019}, \cite{MS_2017}, \cite{SJ_2020}).  The GRM has also attracted some attention in the minimax testing literature, in particular when the collection $b=(b_k)_{1 \leq k \leq m}$ is assumed to be fully known (see, e.g., \cite{Baraud}, \cite{LLM_2012}). \\

In what follows, for simplicity, we focus our attention to the GRM. We adopt a different point of view concerning the uncertainty on the operator $A$. In particular, we assume that the sequence $b= (b_k)_{k \in \mathbb{N}_m}$ (or the corresponding finite rank operator $A$) is unknown but belongs to a given dictionary $\mathcal{B}$ of finite cardinality (i.e., $| \mathcal{B}|<+\infty$). Such a setting is motivated by research investigations in semi-blind astronomical image deconvolution, where one wants to learn both the unknown signal and the unknown convolution kernel (also referred to as the point spread function) from the data. In such a setting, it is usually assumed that the unknown convolution kernel can be well approximated by a linear combination of convolution kernels in a given dictionary of finite cardinality (see, e.g., \cite{FTT_2019}, \cite{TTWBTS_2014}). Our aim, however, is slightly different. Given observations $y=(y_k)_{k \in \mathbb{N}_m}$ from the GRM, where the collection $b= (b_k)_{k \in \mathbb{N}_m}$ is unknown but assumed to belong to a given dictionary $\mathcal{B}$ with $|\mathcal{B}|<+\infty$, we want to address the following classical testing problem
\begin{equation}
H_0: \theta = \theta^0 \quad \mathrm{against} \quad  H_1: \theta \in \Theta,
\label{eq:pbtest}
\end{equation}
where
$$\Theta \subset l^2(N_m)=\left\lbrace \theta \in \mathbb{R}^m: \ \|\theta\|^2 = \sum_{k=1}^m \theta_k^2 < + \infty \right\rbrace$$
such that $\theta^0 \not \in \Theta$. 
In the following, a (non-randomized) test $\Delta:=\Delta(y)$ will be defined as a measurable function of the observations $y=(y_k)_{k \in \mathbb{N}_m}$ from the GRM having values in the set $\lbrace 0,1 \rbrace$. By convention, $H_0$ is rejected if $\Delta=1$ and  $H_0$ is not rejected if $\Delta=0$. Then, given a test $\Delta$, we can investigate
\begin{itemize}
\item The Type I error probability defined as
\begin{equation}
\max_{b\in \mathcal{B}} \mathbb{P}_{\theta^0,b}( \Delta=1),
\label{eq:type1-F}
\end{equation}
which measures the probability to reject $H_0$ when $H_0$ is true (i.e., $\theta=\theta^0$, $b\in \mathcal{B}$). It is often constrained as being bounded by a prescribed level $\alpha \in ]0,1[$. In such a case, $\Delta=\Delta_\alpha$ is called an $\alpha$-level test. 
\item The Type II error probability over a given set $\Theta \subset \ell^2(\mathbb{N}_m)$ defined as
\begin{equation}
\max_{b\in \mathcal{B}}\sup_{\theta\in \Theta} \mathbb{P}_{\theta,b}(\Delta=0),
\label{eq:type2}
\end{equation}
which measures the worst possible probability not to reject $H_0$ when $H_0$ is not true (i.e., when $\theta \in \Theta$, $b\in \mathcal{B}$); one would like to ensure that it is (asymptotically) bounded by a prescribed level $\beta \in ]0,1[$.
\end{itemize}

Note that in the classical setting (i.e., when $b$ is known or partially unknown), the protection against all possible $b \in \mathcal{B}$ (i.e. $\max_{b \in \mathcal{B}}$) is not required. However, in order to produce some kind of robustness when $b$ is not known, we have adapted the definition of Type I and Type II error probabilities to accommodate the (possible) uncertainty on $b \in \mathcal{B}$. \\

Our aim is to provide a description of the smallest possible set $\Theta$ (in a sense which is made precise later on) for which a non-asymptotic control of the Type II error probability is possible. The fact that the collection $b$ is unknown creates some issues. Indeed, any testing strategy is more or less based on a comparison between the observations $y$ and the target signal $\theta^0$. Since we are dealing with an inverse problem, the observations $y$ are not directly connected to the target signal of interest $\theta^0$. In the situation where the `true'  $b$ is unknown, we need to compare $y$ to $\bar b\theta^0$ for any $\bar b\in \mathcal{B}$. In the specific case where $\theta^0=0$ (or is small enough), we have one single image $\bar b \theta^0=0$ (or different images with small enough values) for every $\bar b\in \mathcal{B}$. On the other hand, namely when $\theta^0$ is large enough (in a sense which will be made precise later on), we have a different image $\bar b \theta^0$ for \textit{each} $\bar b \in \mathcal{B}$ (i.e., we have different images with large enough values). In some sense, we have in such a case a composite null hypothesis $H_0$. \\

In this paper, we mainly focus our attention to the specific case where $|\mathcal{B}|=2$. In other words, we will deal with the case where the dictionary contains only two members, one of them is the true one (which gave rise to the data given by model (\ref{eq:gsm})) but we do not know which one. As we will demonstrate, even this simple case investigated in Section \ref{s:minimax}  is quite intricate. In particular, we prove in Section \ref{s:homogeneous} that the separation rates remain the same provided the family $\mathcal{B}$ is homogeneous regarding the null hypothesis $H_0$ (in a sense which is made precise in equation (\ref{eq:cond_reg}) below). On the other hand, the separation rates will deteriorate provided this condition is not satisfied (see Section \ref{s:heterogeneous}). This reveals an interesting phase transition phenomenon regarding the sets $\Theta$ under the alternative. We provide, in both cases, upper and lower bounds. The case $|\mathcal{B}| > 2$, is even more involved and requires different strategies. It is briefly discussed in Section \ref{s:general}. 



\section{Minimax bounds for the case $|\mathcal{B}| = 2$}
\label{s:minimax}

\subsection{The homogeneous regime}
\label{s:homogeneous}

We investigate below  the situation where the dictionary $\mathcal{B}$ is homogeneous, in a sense which is made precise in the Assumption A1 below. To this end, let
$$ \rho_{\epsilon,b} = \epsilon^2 \sqrt{\sum_{j=1}^m b_j^{-4}} \quad \forall b\in \mathcal{B}.$$
We recall that these terms correspond, up to a constant, to the separation rates for a fixed sequence $b$ (see, e.g., \cite{ISS_2012} or \cite{LLM_2012}).\\
\\
\textbf{Assumption A1.} \textit{The set $\mathcal{B}$ and the sequence $\theta^0$ satisfy}
\begin{equation}
\sum_{k=1}^m (b_{k}^{(1)})^{-2} \left(b_k^{(1)} - b_k^{(2)}\right)^2 (\theta_k^{0})^2 \leq  \rho^2_{\epsilon,b^{(1)}} \quad  \forall b^{(1)},b^{(2)} \in \mathcal{B}.
\label{eq:cond_reg}
\end{equation}

The condition (\ref{eq:cond_reg}) indicates, in some sense, that the reconstruction of $\theta^0$ using any member of $\mathcal{B}$ will not be too different compared to the classical detection rate. This is the reason why this setting is called homogeneous. It is trivially satisfied as soon as the dictionary $\mathcal{B}$ contains a single member, i.e., $|\mathcal{B}|=1$. We stress that (\ref{eq:cond_reg}) heavily depends on $\theta_0$, the signal considered under the null hypothesis. Indeed, Assumption A1 will be satisfied by any family $\mathcal{B}$ in the signal detection case, namely when $\theta^0=0$. On the other hand, this condition appears to be quite restrictive when the norm of $\theta^0$ increases, i.e., in the goodness-of-fit setting. \\

Since we have uncertainty on the operator involved, we propose below a testing approach inspired by adaptive methods developed in a non-parametric testing context (see, e.g., \cite{Baraud}, \cite{IngsterI}, \cite{IngsterII}, \cite{IngsterIII}).  Let $\alpha \in ]0,1[$ be fixed.  For any $\bar b\in \mathcal{B}$, we can consider the test $\Delta_{\alpha,\bar b}$ defined as
\begin{equation}
\Delta_{\alpha,\bar b} = \mathbf{1}_{\lbrace \sum_{k=1}^{m} (\bar b_k^{-1} y_k - \theta_k^0)^2 > t_{\alpha,\bar b} \rbrace},
\label{eq:test1}
\end{equation}
where
\begin{equation}
t_{\alpha,\bar b}= \epsilon^2 \sum_{k=1}^{m} \bar b_k^{-2} + c_\alpha \epsilon^2 \sqrt{\sum_{k=1}^{m} \bar b_k^{-4}},
\label{eq:threshold1}
\end{equation}
and $c_\alpha = 1+2 (2\sqrt{x_\alpha} + x_\alpha)$ with $x_\alpha = \ln(1/\alpha)$ for all $\alpha \in ]0,1[$. Note that the procedure defined by (\ref{eq:test1}) and (\ref{eq:threshold1}) is quite standard in the case where the operator is known (see, e.g., \cite{LLM_2011},\cite{LLM_2012}). The underlying idea is to construct an estimator of  $d(\theta,\theta^0)=\| \theta - \theta^0 \|^2$ which will be compared to an appropriate threshold (\ref{eq:threshold1}) that will allow a sharp control for the Type I error probability. In our context, we design a test for each possible candidate $\bar b \in \mathcal{B}$.  \\

Since the true sequence $b \in \mathcal{B}$ is unknown, we can adopt a classical adaptation strategy allowing to select an appropriate candidate from the dictionary $\mathcal{B}$. This strategy consists of aggregating the tests $(\Delta_{\alpha,\bar{b}})_{\bar{b} \in \mathcal{B}}$ defining the following Bonferroni type procedure
\begin{equation}
\Delta^{(1)}_\alpha := \max_{\bar b\in \mathcal{B}} \Delta_{\alpha/2,\bar b}.
\label{eq:D1}
\end{equation}
Such an approach is quite standard in the case where we are faced to an uncertainty in the alternative (e.g., smoothness of the target, underlying dimension of the signal). We point out that our setting is new since the uncertainty on the operator has an impact on both the null and the alternative hypothesis. \\

To describe the behavior of the testing strategy (\ref{eq:D1}), we introduce, for any $b\in \mathcal{B}$ and for any constant $C>0$, the set $\Theta_{1,b}(C)$ defined as
\begin{equation}
\Theta_{1,b}(C ) = \lbrace \theta\in \ell^2(\mathbb{N}_m): \ \| \theta-\theta^0 \|^2 \geq C\rho_{\epsilon,b}^2 \rbrace. 
\label{eq:Theta1}
\end{equation}
The performances of the test $\Delta^{(1)}_\alpha$ are summarized in the following proposition, whose proof is postponed to Section \ref{s:proofs1}.

\begin{proposition}
\label{prop:up1}
Let $\alpha, \beta \in ]0,1[$ be fixed. Assume observations $y=(y_k)_{k \in \mathbb{N}_m}$ from the GRM and consider the testing problem (\ref{eq:pbtest}). Then, provided 
Assumption \textbf{A1} is satisfied, 
$$ \max_{b\in \mathcal{B}} \mathbb{P}_{\theta^0,b}(\Delta_\alpha^{(1)} = 1) \leq \alpha.$$
Moreover
$$ \max_{b\in \mathcal{B}} \sup_{\theta\in \Theta_{1,b}(C^{(1)}_{\alpha,\beta}) } \mathbb{P}_{\theta,b} (\Delta_\alpha^{(1)}=0)\leq \beta,$$
for some constant $C^{(1)}_{\alpha, \beta}>0$ that can be explicitly calculated.
\end{proposition}

As discussed earlier, under condition (\ref{eq:cond_reg}), the images of the target $\theta^0$ for any of the member of the dictionary $\mathcal{B}$ are not far away (in terms of the value of the radius $\rho_{\epsilon,.}$ introduced above). It allows to simultaneously control the Type I error probability of all the tests involved in the aggregation procedure (\ref{eq:D1}).  The condition (\ref{eq:cond_reg}) is, however, not necessary for the control of the Type II error probability. 

Provided condition (\ref{eq:cond_reg}) is satisfied, Proposition \ref{prop:up1} entails that, for any value of the sequence $b\in \mathcal{B}$, it is possible to construct an $\alpha$-level testing procedure $\Delta_\alpha^{(1)}$ with controlled Type II error probability, as soon as $\|\theta-\theta^0\|^2$ is larger than $\rho_{\epsilon,b}$ (up to a constant term). This separation condition is the same as, e.g., the one displayed in \cite{ISS_2012} or \cite{LLM_2012}. In other words, the goodness-of-fit problem is not affected by the uncertainty on the sequence $b$ provided condition (\ref{eq:cond_reg}) holds.  To entails this discussion, the following results provides a simple lower bound. 

\begin{proposition}
\label{prop:lwb1}
Let $\alpha, \beta \in ]0,1[$ such that $\alpha+\beta<1$. Assume observations $y=(y_k)_{k \in \mathbb{N}_m}$ from the GRM and consider the testing problem (\ref{eq:pbtest}). Then
$$ \inf_{\Delta_\alpha} \max_{b\in \mathcal{B}} \sup_{\theta\in \Theta_{1,b}(c^{(1)}_{\alpha,\beta}) } \mathbb{P}_{\theta,b} (\Delta_\alpha=0)> \beta,$$
where $c_{\alpha,\beta}^{(1)} = (2\ln(1+4(1-\alpha-\beta)^2))^{1/4}>0$, and where the infimum is taken over all possible $\alpha$-level tests $\Delta_\alpha$.
\end{proposition}
\textsc{Proof.} Using simple bounds, we have
\begin{eqnarray*}
\inf_{\Delta_\alpha} \max_{b\in \mathcal{B}} \sup_{\theta\in \Theta_{1,b}(c^{(1)}_{\alpha,\beta}) } \mathbb{P}_{\theta,b} (\Delta_\alpha=0)
& \geq & \inf_{\Delta_\alpha} \frac{1}{|\mathcal{B}|} \sum_{b\in \mathcal{B}} \sup_{\theta\in \Theta_{1,b}(c^{(1)}_{\alpha,\beta} )} \mathbb{P}_{\theta,b} (\Delta_\alpha=0) \\
& \geq &  \frac{1}{|\mathcal{B}|} \sum_{b\in \mathcal{B}} \ \inf_{\Delta_\alpha} \sup_{\theta\in \Theta_{1,b}(c^{(1)}_{\alpha,\beta} )} \mathbb{P}_{\theta,b} (\Delta_\alpha=0) > \beta.
\end{eqnarray*}
Indeed, according to \cite{LLM_2012}, we have 
$$\inf_{\Delta_\alpha} \sup_{\theta\in\Theta_{1,b}(c^{(1)}_{\alpha,\beta} )} \mathbb{P}_{\theta,b} (\Delta_\alpha=0) > \beta \quad \forall b\in \mathcal{B}.$$
This concludes the proof. (As it can be seen from the proof, Proposition \ref{prop:lwb1} holds whatever the size of $\mathcal{B}$ is.) 
\begin{flushright}
$\Box$
\end{flushright}

Note that condition (\ref{eq:cond_reg}) leads to a strong constraint on the dictionary $\mathcal{B}$. In some sense, it says that the null hypothesis is not affected by the uncertainty on the operator. In the next section, we discuss the case where this assumption is not satisfied. Such a case corresponds to a situation where the dictionary is not homogeneous.

\subsection{The non-homogeneous regime} 
\label{s:heterogeneous}

We first prove that the condition (\ref{eq:cond_reg}) is mandatory to obtain an alternative governed only by the distance $\|\theta - \theta_0\|^2$. More precisely, if condition (\ref{eq:cond_reg}) is not satisfied (see (\ref{eq:cond_neg}) below), we demonstrate that a more restrictive separation set under the alternative is required. To this end, we introduce, for any $b\in \mathcal{B}$ and for any constant $C>0$, the set $\Theta_{2,b}(C)$ defined as 
\begin{eqnarray*}
\Theta_{2,b} (C)
& = & \left\lbrace \theta\in \ell^2(\mathbb{N}_m):  \| \tilde b^{-1}  (b\theta  - \tilde b \theta^0)\|^2  \geq C \rho_{\epsilon, \tilde b}^2 \quad \forall \tilde b \in \mathcal{B}\setminus \lbrace b \rbrace  \right\rbrace,
\end{eqnarray*}
where 
\[
 \| \tilde b^{-1}  (b\theta  - \tilde b \theta^0)\|^2 = \sum_{k=1}^m \tilde b_k^{-2} (b_k \theta - \tilde b_k \theta^0)^2 \quad \forall b,\tilde b \in \mathcal{B} \quad \mathrm{and} \quad \forall \theta \in \ell^2(\mathbb{N}_m).
\]
Note that if a given $\theta$ belongs to the set $\Theta_{2,b} (C)$, it means that $b\theta$ and $\tilde b \theta^0$ are well separated. This can be understood, in some sense, as an identifiability condition. \\

The following proposition provides a lower bound. The proof is  postponed to Section \ref{s:prop2}.

\begin{proposition} 
\label{prop:2}
Let $\alpha, \beta\in ]0,1[$ be fixed such that $\alpha+\beta <1$. Assume observations $y=(y_k)_{k \in \mathbb{N}_m}$ from the GRM and consider the testing problem (\ref{eq:pbtest}). Assume also that there exist $b^{(1)}, b^{(2)} \in \mathcal{B}$ such that
\begin{equation}
\sum_{k=1}^m (b_{k}^{(1)})^{-2} \left(b_k^{(1)} - b_k^{(2)}\right)^2 (\theta_k^{0})^2 > C_2 \rho^2_{\epsilon,b^{(1)}} 
\label{eq:cond_neg}
\end{equation}
for some constant $C_2\geq 1$. Then,
$$ \inf_{\bar{\Delta}_\alpha} \max_{b\in \mathcal{B}} \sup_{\theta\in \Theta_{1,b}(C_{\alpha,\beta}^{(1)}) \cap \Theta_{2,b}\left(c^{(2)}_{\alpha,\beta}\right)} \mathbb{P}_{\theta,b} (\bar{\Delta}_\alpha=0)> \beta,$$
for some constant $c^{(2)}_{\alpha,\beta}>0$ that can be explicitly calculated, and where the infimum is taken over all possible $\alpha$-level tests $\bar{\Delta}_\alpha$. 
\end{proposition}

Remark first that $\Theta_{1,b}(C_{\alpha,\beta}^{(1)}) \cap \Theta_{2,b}\left(c^{(2)}_{\alpha,\beta}\right) \subset \Theta_{1,b}(C_{\alpha,\beta}^{(1)})$. According to Proposition \ref{prop:2}, we hence obtain
\[
\inf_{\bar{\Delta}_\alpha} \max_{b\in \mathcal{B}} \sup_{\theta\in \Theta_{1,b}(C_{\alpha,\beta}^{(1)}) } \mathbb{P}_{\theta,b} (\bar{\Delta}_\alpha=0)> \beta
\]
Despite we are working on $\Theta_{1,b}(C_{\alpha,\beta}^{(1)})$ (which is known to ensure the separability of $H_0$ and $H_1$ when $b$ is known), every test will hence fail in this specific situation and we have to introduce additional separation conditions.\\ 


We prove below that these associated separation conditions are optimal (up to constants). To this end, we will use a different aggregation strategy. Recall that we cannot use a Bonferroni approach since the condition \ref{eq:cond_reg} (which allows for a uniform control of the Type I error probability) is no longer satisfied. Instead, we use a minimum type aggregation procedure. More formally, let $\Delta^{(2)}_\alpha$ be the testing procedure defined as
$$ \Delta^{(2)}_\alpha := \min_{\bar b\in \mathcal{B}} \Delta_{\alpha,\bar b}, $$
where the tests $\Delta_{\alpha,\bar b}$ have been introduced in (\ref{eq:test1}) and (\ref{eq:threshold1}). It appears that $\Delta^{(2)}_\alpha$ is indeed an $\alpha$-level testing procedure. Indeed, 
\begin{eqnarray}
\max_{b \in \mathcal{B}} \mathbb{P}_{\theta^0,b}(\Delta^{(2)}_{\alpha} = 1)
& = & \max_{b \in \mathcal{B}} \mathbb{P}_{\theta^0,b} \left( \min_{\bar b\in \mathcal{B}} \Delta_{\alpha,\bar b} = 1   \right) \nonumber \\
& = & \max_{b \in \mathcal{B}} \mathbb{P}_{\theta^0,b} \left( \bigcap_{\bar b\in \mathcal{B}} \lbrace \Delta_{\alpha,\bar b} = 1 \rbrace \right) \nonumber \\
& \leq & \max_{b \in \mathcal{B}} \mathbb{P}_{\theta^0,b} (\Delta_{\alpha,b} =1) \leq \alpha,
\label{eq:numproof}
\end{eqnarray}
where, for the last inequality, we used Lemma 2 in \cite{LLM_2012} and standard results when $b \in \mathcal{B}$ is known. Then, a control of the Type II error probability is displayed in the following proposition, whose proof is deferred to Section \ref{s:proof_up_2}.

\begin{proposition}
\label{prop:upper_2}
Let $\alpha, \beta \in [0,1]$ be fixed. Assume observations $y=(y_k)_{k \in \mathbb{N}_m}$ from the GRM and consider the testing problem (\ref{eq:pbtest}). Then,
$$ \max_{b\in \mathcal{B}}  \sup_{\theta\in \Theta_{1,b}(C_{\alpha,\beta}^{(1)}) \cap \Theta_{2,b}\left(C^{(2)}_{\alpha,\beta}\right)} \mathbb{P}_{\theta,b} (\Delta_\alpha^{(2)}=0)\leq \beta,$$
for some constant $C^{(2)}_{\alpha,\beta}>0$ that can be explicitly calculated.
\end{proposition}

We stress that this results holds without any additional restriction. More precisely, it is not necessary that condition (\ref{eq:cond_neg}) holds to get this upper bound. 

Proposition \ref{prop:2} together with Proposition \ref{prop:upper_2} provides the minimax separation condition in this non-homogeneous regime, up to constants. In particular, additional restriction are necessary (namely $\theta \in \Theta_{2,.}$) to control the Type II error probability, provided \eqref{eq:cond_neg} holds.

\subsection{Discussion}
The previous sections provides separation conditions for our testing problem in both homogeneous and non-homogeneous regimes. These regimes are characterised by the value of the term
\[
D_{\theta_0}(b^{(1)},b^{(2)}) = \sum_{k=1}^m (b_{k}^{(1)})^{-2} \left(b_k^{(1)} - b_k^{(2)}\right)^2 (\theta_k^{0})^2 \quad \forall b^{(1)},b^{(2)} \in \mathcal{B},
\]
that corresponds in some sense to some divergence between the sequences involved in $\mathcal{B}$. Both regimes differ from the value of the constant $C_2$ (see \eqref{eq:cond_reg} and \eqref{eq:cond_neg}).  We provide below a short informal discussion concerning the impact of this constant on the separation rates. \\

We can first notice that materials displayed in Section \ref{s:homogeneous} can be adapted to include the divergence terms $D_{\theta_0}(b^{(1)},b^{(2)}) $. In particular, modifying the construction of $\Delta_\alpha^{(1)}$ (by changing the value of the threshold, one may obtain a control of the Type II error provided
\[
\| \theta - \theta_0 \|^2 \geq C_{\alpha,\beta} \rho_{\epsilon,b}^2 + D_{\theta_0}(b,\bar b),
\]
for some constant term $C_{\alpha,\beta} $ and where $b$ denotes the 'true' sequence involved in the model. Provided $D_{\theta_0}(b,\bar b) \leq C_2 \rho_{\epsilon,b}$, we re-discover the rates displayed in Proposition \ref{prop:up1}, replacing $C_{\alpha,\beta}^{(1)}$ by $C_{\alpha,\beta} + C_2$. Then, we can notice that the separation conditions dramatically deteriorates as soon as the value of $C_2$ increase. The distinction between both regimes made in this paper should be understood in this direction: the strategy displayed in Section \ref{s:heterogeneous} might be advantageous for large values of $C_2$. \\

The separation condition exhibited in Section \ref{s:heterogeneous} gathers two sets that have different meanings. The sets $\Theta_1$ provides condition on the distance between the target $\theta$ and the null parameter $\theta_0$ and is quite standard in the testing literature. The second set $\Theta_2$ heavily depends on the heterogeneity of the dictionary and can be related to the divergence $D_{\theta_0}$ in some specific situations. For instance, up to some conditions on the sequence $(\bar b_k^{-1} b_k)_{k\in \mathbb{N}_m}$,
\[
\| \bar b^{-1} (b\theta - \bar b \theta^0\|^2 \sim \| \bar b^{-1} b (\theta-\theta_0) \|^2 + \| (b-\bar b)\theta^0\|^2 \sim \| \theta - \theta^0\|^2 + D_{\theta_0}(\bar b, b),
\]
provided $\| \theta - \theta^0\|^2$ is negligible compared to $D_{\theta_0}(\bar b, b)$. In this latter case, the separation condition associated to the set $\Theta_2$ is mainly driven by the $D_{\theta_0}(\bar b, b)$. In particular, the setting where $D_{\theta_0}(\bar b, b) > C_2 \rho_{\epsilon,\bar b}^2$ for large values of $C_2$ should be handled  via the testing strategy displayed in Section \ref{s:heterogeneous}.

\section{Extension to the case $|\mathcal{B}|>2$}
\label{s:general}

In this section, we briefly discuss the case where $|\mathcal{B}|>2$. We assume that the dictionary can be written as $\mathcal{B} = \mathcal{B}_H \cup \mathcal{B}_{H^c}$, where $\mathcal{B}_H \cap \mathcal{B}_{H^c} = \lbrace \varnothing \rbrace$  and $\mathcal{B}_H$ is such that 
\begin{equation}
\sum_{k=1}^m (b_{k}^{(1)})^{-2} \left(b_k^{(1)} - b_k^{(2)}\right)^2 (\theta_k^{0})^2 \leq  \rho^2_{\epsilon,b^{(1)}} \quad  \forall b^{(1)},b^{(2)} \in \mathcal{B}_H.
\label{eq:cond_reg_gen}
\end{equation}
In other words, we ask for the dictionary $\mathcal{B}$ to be decomposed as the union of a homogeneous part (associated to $\mathcal{B}_H$) and a non-homogeneous part (associated to $\mathcal{B}_H^c$). All the members of $\mathcal{B}_H$ satisfy condition (\ref{eq:cond_reg_gen}) which, in a sense, corresponds to Assumption A1 in the case where $|\mathcal{B}|=2$. \\

The above decomposition and the results for the case $|\mathcal{B}|=2$ help us to construct a testing procedure.  First, define
$$ \Delta_{\alpha,\mathcal{B}_H} = \max_{b\in \mathcal{B}_H} \Delta_{\frac{\alpha}{|\mathcal{B}_H|}, b},$$
and then consider
\begin{equation} 
\Delta_\alpha = \min \left( \Delta_{\alpha,\mathcal{B}_H} , \min_{b\in \mathcal{B}_{H^c}} \Delta_{\alpha,b}\right).
\label{eq:general_case_test}
\end{equation}
In some sense, the test $\Delta_\alpha$ introduced in (\ref{eq:general_case_test}) corresponds to the aggregation of $ \Delta_{\alpha,\mathcal{B}_H}$ (which is the analogue of $\Delta_\alpha^{(1)}$ in Section \ref{s:homogeneous}) and of $ \min_{b\in \mathcal{B}_{H^c}} \Delta_{\alpha,b}$ (which is the analogue of $\Delta_\alpha^{(2)}$ in Section \ref{s:heterogeneous}).\\

Now, we introduce the following separation sets. For any $b,\bar b \in \mathcal{B}$, and for any constant $C>0$, define
$$ \Theta_{3,b}^{\bar b}(C) = 
 \left\lbrace \theta\in \ell^2(\mathbb{N}_m):  \| \bar b^{-1}  (b\theta  - \bar b \theta^0)\|^2  \geq C\rho_{\epsilon,\bar b}^2 \right\rbrace,$$
and consider, for any constants $C_1>0$, $C_2>0$, the set
$$  \Theta_{3,b} (C_1,C_2)= \bigcup_{\bar b \in \mathcal{B}_H}  \Theta_{3,b}^{\bar b}(C_1) \cap  \bigcap_{\bar b \in \mathcal{B}_H^c}  \Theta_{3,b}^{\bar b}(C_2).$$
The performance of the test $\Delta_\alpha$ is summarized in the following proposition, whose proof is postponed to Section \ref{s:general_B}.

\begin{proposition}
\label{prop:up4}
Let $\alpha, \beta \in ]0,1[$ be fixed. Assume observations $y=(y_k)_{k \in \mathbb{N}_m}$ from the GRM and consider the testing problem (\ref{eq:pbtest}). Then, provided 
condition (\ref{eq:cond_reg_gen}) is satisfied, we have
$$ \max_{b\in \mathcal{B}} \mathbb{P}_{\theta^0,b}(\Delta_\alpha = 1) \leq \alpha.$$
Moreover
$$ \max_{b\in \mathcal{B}} \sup_{\theta\in \Theta_{3,b}\big(C_{\frac{\alpha}{|\mathcal{B}_H|},\beta},C_{\alpha,\frac{\beta}{|\mathcal{B}_H^c|}}\big)} \mathbb{P}_{\theta,b} (\Delta_\alpha=0)\leq \beta,$$
for some constants $C_{\frac{\alpha}{|\mathcal{B}_H|},\beta}>0$ and $C_{\alpha,\frac{\beta}{|\mathcal{B}_H^c|}}>0$ that can be explicitly calculated.
\end{proposition}

The separation set $\Theta_{3,b}$ provides strong separation conditions. In particular, the larger the $\mathcal{B}_H^c$, the smaller the $\Theta_{3,b}$. In other words, if the non-homogeneous part $\mathcal{B}_H^c$ of the dictionary $\mathcal{B}$ is too large, then, the goodness-of-fit problem is quite complicated and this affects the separation rates. \\

We point out that the decomposition of the dictionary as $\mathcal{B} = \mathcal{B}_H \cup \mathcal{B}_{H^c}$ is not unique. Indeed, $\mathcal{B}_H$ gathers sequences $b\in \mathcal{B}$ that are close to each other (in the sense of condition (\ref{eq:cond_reg_gen})). Even in the case where $|\mathcal{B}|=3$, constructing $\mathcal{B}_H$ is not trivial. Moreover, remark that the separation set $\Theta_{3,b}\big(C_{\frac{\alpha}{|\mathcal{B}_H|},\beta},C_{\alpha,\frac{\beta}{|\mathcal{B}_H^c|}}\big)$ depends explicitly on this construction. The choice of $\mathcal{B}_H$ might be driven by the goal of maximizing the size of $\Theta_{3,b}\big(C_{\frac{\alpha}{|\mathcal{B}_H|},\beta},C_{\alpha,\frac{\beta}{|\mathcal{B}_H^c|}}\big)$ (hence providing the weakest possible separation condition).  We also mention that we could use a classical adaptation strategy to choose a data driven set $\mathcal{B}_H \subset \mathcal{B}$. The idea is to associate to each $\mathcal{B}_H$ a corresponding test $\Delta_\alpha = \Delta_\alpha(\mathcal{B}_H)$ and then use the classical Bonferonni aggregation procedure $\Delta_\alpha^*$ defined as
$$ \Delta_\alpha^* = \max_{\mathcal{B}_H \in \mathfrak{B}} \Delta_{\frac{\alpha}{|\mathfrak{B}|}}(\mathcal{B}_H),$$
where $\mathfrak{B}$ denotes all the subsets of $\mathcal{B}$ that satisfies condition (\ref{eq:cond_reg_gen}). \\

In this context, deriving the optimality of the separation set $\Theta_{3,b}\big(C_{\frac{\alpha}{|\mathcal{B}_H|},\beta},C_{\alpha,\frac{\beta}{|\mathcal{B}_H^c|}}\big)$ in Proposition \ref{prop:up4} is quite intricate and requires to design a lower bound. We hope to address it in a future work.

\section{Proofs and technical results}
\label{s:proofs}

\subsection{Proof of Proposition \ref{prop:up1}}
\label{s:proofs1}

Below, we investigate both Type I and Type II error probabilities.  We first concentrate our attention on the Type I error probability. Using simple bounds,
\begin{eqnarray*}
\max_{b\in \mathcal{B}} \mathbb{P}_{\theta^0,b} (\Delta^{(1)}_\alpha = 1)
& = & \max_{b\in \mathcal{B}} \mathbb{P}_{\theta^0,b} \left( \max_{\bar b\in \mathcal{B}} \Delta_{\alpha/2,\bar b}=1\right)\\
& = & \max_{b\in \mathcal{B}} \mathbb{P}_{\theta^0,b} \left( \bigcup_{\bar b\in \mathcal{B}} \lbrace \Delta_{\alpha/2,\bar b}=1 \rbrace\right) \\
& \leq &\max_{b\in \mathcal{B}}  \sum_{\bar b\in \mathcal{B}} \mathbb{P}_{\theta^0,b} \left( \Delta_{\alpha/2,\bar b}=1 \right) \\
& \leq & \alpha,
\end{eqnarray*}
provided
\begin{equation}
\mathbb{P}_{\theta^0,b} \left( \Delta_{\alpha/2,\bar b}=1 \right) \leq \frac{\alpha}{2} \quad \forall\, b,\bar b\in \mathcal{B}.
\label{eq:toprove}
\end{equation}
Lemma \ref{lem:1} below provides conditions under which (\ref{eq:toprove}) holds. We can now state an upper bound on the Type II error probability. We simply remark that for any $\theta \not = \theta^0$ and for any $b\in \mathcal{B}$,
$$  \mathbb{P}_{\theta,b} (\Delta^{(1)}_{\alpha} =0) = \mathbb{P}_{\theta,b} \left( \max_{\bar b\in \mathcal{B}} \Delta_{\alpha/2,\bar b}=0\right)
= \mathbb{P}_{\theta,b} \left( \bigcap_{\bar b\in \mathcal{B}} \lbrace \Delta_{\alpha/2,\bar b}=0 \rbrace\right) \leq \mathbb{P}_{\theta,b} (\Delta_{\alpha/2,b}=0) . $$
In other words, the  Type II error probability of the aggregated test $\Delta^{(1)}_\alpha$ is bounded from above by the Type II error probability of the standard test when $b \in \mathcal{B}$ is known.  We can hence apply Proposition 2 in \cite{LLM_2012} to get
$$ \max_{b\in \mathcal{B}} \sup_{\theta\in \Theta_{1,b}(C^{(1)}_{\alpha,\beta} )} \mathbb{P}_{\theta,b} (\Delta_\alpha^{(1)}=0)\leq  \max_{b\in \mathcal{B}} \sup_{\theta\in \Theta_{1,b}(C^{(1)}_{\alpha,\beta} )} \mathbb{P}_{\theta,b} (\Delta_{\alpha/2,b}=0) \leq \beta,$$
where
\begin{equation}
C^{(1)}_{\alpha,\beta}=\sqrt{2x_{\beta}}+\sqrt{2(x_{\alpha/2}+x_{\beta}) + \sqrt{2}(\sqrt{x_{\alpha/2}}+\sqrt{x_\beta})^{1/2}},
\label{eq:constantC1}
\end{equation}
with $x_\gamma =\ln(1/\gamma)$ for any $\gamma \in ]0,1[$.  This completes the proof of the proposition.

\begin{flushright}
$\Box$
\end{flushright}

\begin{lemma}
\label{lem:1}
Assume observations $y=(y_k)_{k \in \mathbb{N}_m}$ from the GRM and consider the testing problem (\ref{eq:pbtest}). Assume that $\theta^0 \in \ell^2(\mathbb{N}_m)$ is such that condition (\ref{eq:cond_reg}) is satisfied.  Then, (\ref{eq:toprove}) holds.
\end{lemma}
\textsc{Proof}. Let $b,\bar b\in \mathcal{B}$ be fixed. Then
$$ \mathbb{P}_{\theta^0,b} \left( \Delta_{\alpha/2,\bar b}=1 \right)  = \mathbb{P}_{\theta^0,b} \left( \sum_{k=1}^m \left(\bar b_k^{-1} y_k - \theta_k^0\right)^2 > t_{\alpha/2,\bar b}   \right) = \mathbb{P} ( T_{\bar b}^0  > t_{\alpha/2,\bar b} ),$$
with
$$ T_{\bar b}^0 := \sum_{k=1}^m \left(\bar b_k^{-1} y_k - \theta_k^0\right)^2 = \sum_{k=1}^m \left( \left(\frac{b_k}{\bar b_k}-1\right) \theta_k^0 + \epsilon \bar b_k^{-1} \xi_k \right)^2.$$
Note that
$$ \mathbb{E}[T_{\bar b}^0] =  \sum_{k=1}^m \left(\frac{b_k}{\bar b_k}-1\right)^2 (\theta_k^0)^2 + \epsilon^2\sum_{k=1}^m \bar b_k^{-2}.$$
Following Lemma 2 in \cite{LLM_2012}, we get
$$ \mathbb{P}\left( T_{\bar b}^0 -\mathbb{E}[T_{\bar b}^0]  > 2\sqrt{x_{\alpha/2}} \sqrt{ \epsilon^4 \sum_{k=1}^m \bar b_k^{-4} + 2\epsilon^2 \sum_{k=1}^m \bar b_k^{-2} \left(\frac{b_k}{\bar b_k}-1\right)^2 (\theta_k^0)^2}+ 2 x_{\alpha/2} \epsilon^2 \max_{1 \leq k\leq m} (\bar b_k^{-2}) \right) \leq \frac{\alpha}{2}.$$
Hence,
$$  \mathbb{P}_{\theta^0,b} \left( \Delta_{\alpha/2,\bar b}=1 \right) = \mathbb{P} ( T_{\bar b}^0 -\mathbb{E}[T_{\bar b}^0] > t_{\alpha/2,\bar b}- \mathbb{E}[T_{\bar b}^0] ) \leq \frac{\alpha}{2} $$
as soon as
$$ t_{\alpha/2,\bar b}- \mathbb{E}[T_{\bar b}^0] \geq 2\sqrt{x_{\alpha/2}} \sqrt{ \epsilon^4 \sum_{k=1}^m \bar b_k^{-4} + 2\epsilon^2 \sum_{k=1}^m \bar b_k^{-2} \left(\frac{b_k}{\bar b_k}-1\right)^2 (\theta_k^0)^2}+ 2 \epsilon^2 \max_{1 \leq k\leq m} (\bar b_k^{-2}) x_{\alpha/2}$$
or, equivalently, as soon as
\begin{eqnarray*}
\sum_{k=1}^m \left(\frac{b_k}{\bar b_k}-1\right)^2 (\theta_k^0)^2 &+& \epsilon^2\sum_{k=1}^m \bar b_k^{-2} +  2\sqrt{x_{\alpha/2}} \sqrt{ \epsilon^4 \sum_{k=1}^m \bar b_k^{-4} + 2\epsilon^2 \sum_{k=1}^m \bar b_k^{-2} \left(\frac{b_k}{\bar b_k}-1\right)^2 (\theta_k^0)^2} \\
&+& 2 \epsilon^2x_{\alpha/2} \max_{1 \leq k\leq m} (\bar b_k^{-2}) \leq t_{\alpha/2,\bar b}.
\end{eqnarray*}
Applying the standard inequalities $\sqrt{a+b} \leq \sqrt{a} + \sqrt{b}$ for any $a,b \in \mathbb{R_+}$ and $2ab \leq a^2 + b^2$ for any $a,b \in \mathbb{R}$, and since
$ \max_{1 \leq k\leq m} (\bar b_k^{-2}) \leq (\sum_{k=1}^m \bar b_k^{-4})^{1/2}$,
we can see that this condition holds provided
$$ (1+ 2\sqrt{x_{\alpha/2}})\sum_{k=1}^m \left(\frac{b_k}{\bar b_k}-1\right)^2 (\theta_k^0)^2 + \epsilon^2\sum_{k=1}^m \bar b_k^{-2} +  2(\sqrt{x_{\alpha/2}}+ x_{\alpha/2}) \epsilon^2 \sqrt{\sum_{k=1}^m \bar b_k^{-4}} \leq t_{\alpha/2,\bar b}.$$
Setting
$$ t_{\alpha/2,\bar b} = \epsilon^2 \sum_{k=1}^m \bar b_k^{-2}  + (1+  2(2\sqrt{x_{\alpha/2}}+ x_{\alpha/2}) )\epsilon^2 \sqrt{\sum_{k=1}^m \bar b_k^{-4}} ,$$
or, equivalently,
$$ t_{\alpha/2,\bar b} = \epsilon^2 \sum_{k=1}^m \bar b_k^{-2}  + c_{\alpha/2} \epsilon^2 \sqrt{\sum_{k=1}^m \bar b_k^{-4}} ,
$$
we get that (\ref{eq:toprove}) holds as soon as (\ref{eq:cond_reg}) holds. This completes the proof of Lemma \ref{lem:1}.

\begin{flushright}
$\Box$
\end{flushright}

\subsection{Proof of Proposition \ref{prop:2}}
\label{s:prop2}

\subsubsection{A generic result}
Let $\Theta \subset l^2(\mathbb{N}_m)$. In the following, $\Delta_\alpha$ denotes a generic $\alpha$-level test satisfying in particular
\begin{equation}
\max_{b \in \mathcal{B}} \mathbb{P}_{\theta^0,b}(\Delta_\alpha = 1) \leq \alpha.
\label{eq:T1lb}
\end{equation}
Then,  with an abuse of notation, for any given $b \in \mathcal{B}$,
\begin{eqnarray*}
\inf_{\Delta_\alpha} \max_{b \in \mathcal{B}} \sup_{\theta \in \Theta} \mathbb{P}_{\theta,b}(\Delta_\alpha =0)
& \geq & \inf_{\Delta_\alpha}  \sup_{\theta \in \Theta} \mathbb{P}_{\theta,b}(\Delta_\alpha =0) \\
& \geq & \inf_{\Delta_\alpha}  \int_{\Theta} \mathbb{P}_{\theta,b}(\Delta_\alpha =0)d\mu(\theta)\\
& = & \inf_{\Delta_\alpha} \mathbb{P}_{\mu,b}(\Delta_\alpha =0),
\end{eqnarray*}
where $\mu$ is a probability measure on $\Theta$ and
$ \mathbb{P}_{\mu,b} = \int_\Theta \mathbb{P}_{\theta,b} d\mu(\theta)$.
Then, for any $\bar b\in \mathcal{B}$, using (\ref{eq:T1lb}),
\begin{eqnarray*}
\inf_{\Delta_\alpha} \max_{b \in \mathcal{B}} \sup_{\theta \in \Theta} \mathbb{P}_{\theta,b}(\Delta_\alpha =0)
& \geq & \inf_{\Delta_\alpha} \left[   \mathbb{P}_{\theta^0,\bar b}(\Delta_\alpha=0) + \mathbb{P}_{\mu,b}(\Delta_\alpha =0) - \mathbb{P}_{\theta^0,\bar b}(\Delta_\alpha=0) \right]  \\
& \geq & \inf_{\Delta_\alpha} \left[   (1- \alpha) + \mathbb{P}_{\mu,b}(\Delta_\alpha =0) - \mathbb{P}_{\theta^0,\bar b}(\Delta_\alpha=0)  \right] \\
& \geq & \inf_{\Delta_\alpha} \left[   (1- \alpha) - \left| \mathbb{P}_{\mu,b}(\Delta_\alpha =0) - \mathbb{P}_{\theta^0,\bar b}(\Delta_\alpha=0) \right| \right] \\
& \geq & (1-\alpha) - \sup_{A\in \mathcal{A}} \left| \mathbb{P}_{\mu,b}(A) - \mathbb{P}_{\theta^0,\bar b}(A) \right| \\
& = & (1-\alpha) - \frac{1}{2} \| \mathbb{P}_{\mu,b} - \mathbb{P}_{\theta^0,\bar b} \|,
\end{eqnarray*}
where $\| \mathbb{P}_{\mu,b} - \mathbb{P}_{\theta^0,\bar b} \|$ denotes the total variation distance between the probability measures $\mathbb{P}_{\mu,b}$ and $\mathbb{P}_{\theta^0,\bar b}$. Whenever $\mathbb{P}_{\mu,b}$ is absolutely continuous with respect to $\mathbb{P}_{\theta^0,\bar b}$, the total variation norm between these two probabilities can be easily computed. Setting
$$ L_{\mu}(y) = \frac{d\mathbb{P}_{\mu,b}}{d\mathbb{P}_{\theta^0,\bar b}}(y),$$
and following standard arguments (see, eg, \cite{Baraud}, \cite{IS_2003}, \cite{MS_2015}), we get
$$ \| \mathbb{P}_{\mu,b} - \mathbb{P}_{\theta^0,\bar b} \| \leq \left(  \mathbb{E}_{\theta^0,\bar b}[L^2_\mu(y)] - 1   \right)^{1/2}.$$
Hence, for $\mathcal{B} = \lbrace b,\bar b \rbrace$,
$$ \inf_{\Delta_\alpha} \max_{b\in \mathcal{B}} \sup_{\theta \in \Theta} \mathbb{P}_{\theta,b}(\Delta_\alpha =0) \geq 1-\alpha - \frac{1}{2} \left(  \mathbb{E}_{\theta^0,\bar b}[L^2_\mu(y)] - 1   \right)^{1/2}.$$
Our aim is now to make the quantity in the right hand side of the previous inequality larger than $\beta$, which is equivalent to establish that
$$ \mathbb{E}_{\theta^0,\bar b}[L^2_\mu(y)] < C_{\alpha,\beta},$$
where $C_{\alpha,\beta}=1+4(1-\alpha-\beta)^2$. We calculate below the value of $\mathbb{E}_{\theta^0,\bar b} [L^2_{\mu}(y)] $. We have
\begin{eqnarray}
\mathbb{E}_{\theta^0,\bar b} [L^2_{\mu}(y)]
& = & \mathbb{E}_{\theta^0,\bar b} \left[\left(\frac{d\mathbb{P}_{\mu,b}}{d\mathbb{P}_{\theta^0,\bar b}}(y)\right)^2\right]  \nonumber \\
& = & \mathbb{E}_{\theta^0,\bar b} \left[\left( \mathbb{E}_{\mu} \left[ \prod_{k=1}^m \frac{\frac{1}{\sqrt{2\pi\epsilon^2}} \exp \left( - \frac{(y_k - b_k \theta_k)^2}{2\epsilon^2} \right)}{\frac{1}{\sqrt{2\pi\epsilon^2}} \exp \left( - \frac{(y_k - \bar b_k \theta_k^0)^2}{2\epsilon^2} \right)}   \right]  \right)^2\right]    \nonumber \\
& = & \mathbb{E}_{\theta^0,\bar b} \left[\left(   \mathbb{E}_{\mu} \left[ \prod_{k=1}^m  \exp \left(  \frac{1}{2\epsilon^2} \left\lbrace   (y_k - \bar b_k\theta_k^0)^2 - (y_k - b_k\theta_k)^2     \right\rbrace  \right) \right] \right)^2\right]    \nonumber  \\
& = & \mathbb{E} \left[\left(    \mathbb{E}_{\mu} \left[ \prod_{k=1}^m  \exp \left(  \frac{1}{2\epsilon^2} \left\lbrace   \epsilon^2 \xi_k^2 - ((\bar b_k\theta_k^0 - b_k\theta_k) + \epsilon \xi_k)^2     \right\rbrace  \right) \right]  \right)^2\right]     \nonumber   \\
& = & \mathbb{E} \left[\left(  \mathbb{E}_{\mu} \left[ \prod_{k=1}^m  \exp \left(  -\frac{(\bar b_k\theta_k^0 - b_k\theta_k)^2}{2\epsilon^2} - \frac{\xi_k(\bar b_k\theta_k^0 - b_k\theta_k)}{\epsilon}  \right) \right]  \right)^2\right],
\label{eq:LRlb}
\end{eqnarray}
where we have substituted $y_k = \bar b_k \theta_k^0 + \epsilon \xi_k$ to take into account the expectation with respect to $\mathbb{P}_{\theta^0,\bar b}$ as defined above.\\

We propose below a construction for the measure $\mu$ to obtain the lower bound for the set $\Theta = \Theta_{2,b}(c^{(2)}_{\alpha,\beta} \rho_{\epsilon,\mathcal{B}})$ displayed in Proposition \ref{prop:2}.

\subsubsection{Proof of Proposition \ref{prop:2}}
We assume here that $\theta^0 \in \ell^2(\mathbb{N}_m)$ and consider $b, \bar b \in \mathcal{B}$ are such that
\begin{equation}
\sum_{k=1}^m b_{k}^{-2} (b_k - \bar b_k)^2 (\theta_k^{0})^2 > C_2 \rho^2_{\epsilon,b},
\label{eq:cond4}
\end{equation}
for some constant $C_2>0$.  Now, for some $\tau\in [0,1]$ and $R_{\epsilon,\bar b}>0$ (whose values will be defined later on), we define the sequence $\theta=(\theta_k)_{k\in \mathbb{N}_m}$ as
\begin{equation}
b_k \theta_k = \bar b_k \theta_k^0 + \omega_k \frac{\tau R_{\epsilon,\bar b} \bar b_k^{-1}}{\sqrt{\sum_{k=1}^m \bar b_k^{-4}}},  \quad k=1,\dots ,m,
\label{eq:theta0I}
\end{equation}
where the $\omega_k$, $k=1,2,\ldots,m$, are independent Rademacher random variables, i.e., $\mathbb{P}(\omega_k=1)=\mathbb{P}(\omega_k=-1)=1/2$, $k=1,2,\ldots,m$.  Using (\ref{eq:LRlb}), we get
\begin{eqnarray*}
\mathbb{E}_{\theta^0,\bar b} [L^2_{\mu}(y)]
& = & \mathbb{E} \left[\left(  \mathbb{E}_{\omega} \left[ \prod_{k=1}^m  \exp \left(  -\frac{(\omega_k \tau R_{\epsilon,\bar b} \bar b_k^{-1})^2}{2\epsilon^2 \sum_{j=1}^m \bar b_j^{-4}} - \frac{\xi_k\omega_k \tau R_{\epsilon,\bar b} \bar b_k^{-1}}{\epsilon\sqrt{\sum_{j=1}^m \bar b_j^{-4}}}  \right) \right]  \right)^2\right] \\
& = & \mathbb{E} \left[  \exp \left( - \frac{\tau^2 R_{\epsilon,\bar b}^2 \sum_{k=1}^m \bar b_k^{-2}}{\epsilon^2 \sum_{j=1}^m \bar b_j^{-4}}\right) \prod_{k=1}^m \cosh^2 \left( \frac{\xi_k \tau R_{\epsilon,\bar b} \bar b_k^{-1}}{\epsilon\sqrt{\sum_{j=1}^m \bar b_j^{-4}}}  \right) \right]\\
& = & \exp \left( - \frac{\tau^2 R_{\epsilon,\bar b}^2 \sum_{k=1}^m \bar b_k^{-2}}{\epsilon^2 \sum_{j=1}^m \bar b_j^{-4}}\right) \prod_{k=1}^m \mathbb{E} \left[\cosh^2 \left( \frac{\xi_k \tau R_{\epsilon,\bar b} \bar b_k^{-1}}{\epsilon\sqrt{\sum_{j=1}^m \bar b_j^{-4}}}  \right) \right].
\end{eqnarray*}
Using the well known inequalities
$$ \mathbb{E}[ \cosh^2(\lambda Z)] = e^{\lambda^2} \cosh^2(\lambda^2)\quad \mathrm{and} \quad \cosh(x) \leq e^{\frac{x^2}{2}},$$
for any $\lambda,x \in \mathbb{R}$ and $Z \sim \mathcal{N}(0,1)$, we get
\begin{eqnarray*}
\mathbb{E}_{\theta^0,\bar b} [L^2_{\mu}(y)]
& = & \prod_{k=1}^m \cosh^2 \left( \frac{\tau^2 R_{\epsilon,\bar b}^2 \bar b_k^{-2}}{\epsilon^2\sum_{j=1}^m \bar b_j^{-4}} \right) \\
& \leq &  \prod_{k=1}^m \exp \left( \frac{\tau^4 R_{\epsilon,\bar b}^4 \bar b_k^{-4}}{\epsilon^4 \left(\sum_{j=1}^m \bar b_j^{-4}\right)^2} \right) \\
& = & \exp \left( \frac{\tau R_{\epsilon,\bar b}^4}{\epsilon^4 \sum_{j=1}^m \bar b_j^{-4}} \right) <  C_{\alpha,\beta},
\end{eqnarray*}
provided that
\begin{equation}
\tau^2 R_{\epsilon,\bar b}^2  < \sqrt{\ln(C_{\alpha,\beta})} \epsilon^2 \sqrt{\sum_{j=1}^m \bar b_j^{-4}}.
\label{eq:LBIPcond1}
\end{equation}
Taking $R_{\epsilon,\bar b}^2=\sqrt{\ln(C_{\alpha,\beta})}\rho_{\epsilon,\bar b}^2$, the above inequality holds for any $\tau \in [0,1]$.\\ 

To conclude the proof, 
we need to show that, for some $c^{(2)}_{\alpha,\beta}>0$, the sequence $\theta = (\theta_k)_{k\in \mathbb{N}_m}$ constructed in (\ref{eq:theta0I}) belongs to $ \Theta_{1,b}(C_{\alpha,\beta}^{(1)}) \cap \Theta_{2,b}\left(c^{(2)}_{\alpha,\beta}\right)$. In other words, we need to show that
\begin{equation}
\| \theta - \theta^0 \|^2 \geq C^{(1)}_{\alpha,\beta} \rho^2_{\epsilon,b} \quad \text{and} \quad
\| \bar{b}^{-1}(b\theta-\bar{b}\theta^{0})\|^2 \geq c^{(2)}_{\alpha,\beta} \rho_{\epsilon,\bar b}^2,
\label{eq:lwb_requirement}
\end{equation}
where $C^{(1)}_{\alpha,\beta}$ is given in (\ref{eq:constantC1}). 
First note that, by the construction (\ref{eq:theta0I}),
$$ \sum_{k=1}^{m} \bar b_k^{-2}(b_k\theta_k-\bar b_k\theta_k^0)^2 = \tau^2 R_{\epsilon,\bar b}^2 = \tau^2 \sqrt{\ln(C_{\alpha,\beta})}\rho_{\epsilon,\bar b}^2 = c^{(2)}_{\alpha,\beta} \rho_{\epsilon,\bar b}^2,$$
where
\begin{equation}
c^{(2)}_{\alpha,\beta} = \tau^2 \sqrt{\ln(C_{\alpha,\beta})}.
\label{eq:lb_constant}
\end{equation}
Moreover, using the inequality $(a+b)^2 \geq (1-\gamma)a^2 + (1- \gamma^{-1})b^2$ for any $a,b \in \mathbb{R}$ and any $\gamma \in ]0,1[$,
\begin{eqnarray*}
\sum_{k=1}^{m}(\theta_k-\theta_k^{0})^2&=& \sum_{k=1}^m \left( b_k^{-1} (\bar b_k - b_k)\theta_k^0 + \omega_k \frac{\tau R_{\epsilon,\bar b} b_k^{-1}\bar b_k^{-1}}{\sqrt{\sum_{k=1}^m \bar b_k^{-4}}}  \right)^2\\
& \geq & (1-\gamma)  \sum_{k=1}^m b_k^{-2} (\bar b_k - b_k)^2 (\theta_k^0)^2 + (1- \gamma^{-1})  \tau^2 R^2_{\epsilon,\bar b}\frac{\sum_{k=1}^m b_k^{-2}\bar b_k^{-2}}{\sum_{k=1}^m \bar b_k^{-4}} \\
&\geq& (1-\gamma)  \sum_{k=1}^m b_k^{-2} (\bar b_k - b_k)^2 (\theta_k^0)^2
+ (1- \gamma^{-1})  c^{(2)}_{\alpha,\beta} \frac{\epsilon^2 \sum_{k=1}^m b_k^{-2}\bar b_k^{-2}}{\sqrt{\sum_{k=1}^m \bar b_k^{-4}}}
\\
&\geq& (1-\gamma) C_2 \rho^2_{\epsilon,b} + (1- \gamma^{-1}) c^{(2)}_{\alpha,\beta} \epsilon^2 \sqrt{\sum_{k=1}^m b_k^{-4}} \\
&=& \left((1-\gamma) C_2 + (1- \gamma^{-1})\tau^2 \sqrt{\ln(C_{\alpha,\beta})}\right) \rho^2_{\epsilon,b}
\end{eqnarray*}
where, in the last inequality, we have used (\ref{eq:cond4}), the Cauchy-Schwartz inequality, and the fact that $1-\gamma^{-1}<0$ for any $\gamma \in ]0,1[$. 
The first condition of (\ref{eq:lwb_requirement}) holds as soon as 
$$  \left((1-\gamma) C_2 + (1- \gamma^{-1})\tau^2 \sqrt{\ln(C_{\alpha,\beta})}\right)  > C_{\alpha,\beta}^{(1)},$$
which holds true as soon as $C_2$ is large enough and $\tau$ is small enough. This concludes the proof of the proposition.

\begin{flushright}
$\Box$
\end{flushright}

\subsection{Proof of Proposition \ref{prop:upper_2}} 
\label{s:proof_up_2}
For any $\theta \in \ell^2(\mathbb{N}_m)$ and any $b \in \mathcal{B}$, we have
$$ \mathbb{P}_{\theta,b} (\Delta^{(2)}_{\alpha} = 0) = \mathbb{P}_{\theta,b} \left( \min_{\bar b\in \mathcal{B}} \Delta_{\alpha,\bar b} = 0   \right) = \mathbb{P}_{\theta,b} \left( \bigcup_{\bar b\in \mathcal{B}} \lbrace \Delta_{\alpha,\bar b} = 0 \rbrace \right) \leq \sum_{\bar b\in \mathcal{B}} \mathbb{P}_{\theta,b} \left( \Delta_{\alpha,\bar b} = 0  \right) \leq \beta$$
provided that
$$ \mathbb{P}_{\theta,b} \left( \Delta_{\alpha,\bar b} = 0  \right) \leq \frac{\beta}{2} \quad \forall\; b, \bar b\in \mathcal{B}.$$
From now on, let $b, \bar b\in \mathcal{B}$ be fixed. Then
$$ \mathbb{P}_{\theta,b} \left( \Delta_{\alpha,\bar b}= 0  \right) = \mathbb{P}_{\theta,b} \left( \sum_{k=1}^{m} (\bar b_k^{-1} y_k - \theta_k^0)^2 \leq t_{\alpha,\bar b} \right)= \mathbb{P} \left( T_{\bar b} \leq t_{\alpha,\bar b} \right),$$
where
$$ T_{\bar b} =  \sum_{k=1}^{m} \left( \frac{b_k}{\bar b_k} \theta_k -  \theta_k^0 + \epsilon \bar b_k^{-1} \xi_k  \right)^2. $$
Since
$$ \mathbb{E} [ T_{\bar b}] = \sum_{k=1}^{m} \left( \frac{b_k}{\bar b_k} \theta_k - \theta_k^0 \right)^2 + \epsilon^2 \sum_{k=1}^{m} \bar b_k^{-2},$$
we get, using Lemma 2 in \cite{LLM_2012},
$$ \mathbb{P}_{\theta,b} \left( \Delta_{\alpha,\bar b}= 0  \right) = \mathbb{P} \left( T_{\bar b} - \mathbb{E}[T_{\bar b}] \leq t_{\alpha,\bar b} - \sum_{k=1}^{m} \left( \frac{b_k}{\bar b_k} \theta_k - \theta_k^0 \right)^2 - \epsilon^2 \sum_{k=1}^{m} \bar b_k^{-2} \right) \leq \frac{\beta}{2},$$
provided that,
$$
t_{\alpha,\bar b} - \sum_{k=1}^{m} \left( \frac{b_k}{\bar b_k} \theta_k - \theta_k^0 \right)^2 - \epsilon^2 \sum_{k=1}^{m} \bar b_k^{-2} \leq - 2\sqrt{x_{\beta/2}} \sqrt{ \epsilon^4 \sum_{k=1}^{m} \bar b_k^{-4} + 2 \sum_{k=1}^{m} \epsilon^2 \bar b_k^{-2}\left( \frac{b_k}{\bar b_k} \theta_k - \theta_k^0 \right)^2}
$$
or, equivalently,
$$
\sum_{k=1}^{m} \left( \frac{b_k}{\bar b_k} \theta_k - \theta_k^0 \right)^2 \geq t_{\alpha,\bar b} - \epsilon^2 \sum_{k=1}^{m} \bar b_k^{-2} + 2\sqrt{x_{\beta/2}} \sqrt{ \epsilon^4 \sum_{k=1}^{m} \bar b_k^{-4} + 2 \sum_{k=1}^{m} \epsilon^2 \bar b_k^{-2}\left( \frac{b_k}{\bar b_k} \theta_k - \theta_k^0 \right)^2}.
$$
Applying the standard inequalities $\sqrt{a+b} \leq \sqrt{a} + \sqrt{b}$ for any $a,b \in \mathbb{R}_+$ and $ab \leq a^2/(2\delta) + \delta b^2/2$ for any $a,b \in \mathbb{R}$ and any $\delta \in \mathbb{R}_+$, the above inequality holds as soon as
$$
\sum_{k=1}^{m} \left( \frac{b_k}{\bar b_k} \theta_k - \theta_k^0 \right)^2 \geq t_{\alpha,\bar b} - \epsilon^2 \sum_{k=1}^{m} \bar b_k^{-2} + 2\sqrt{x_{\beta/2}} \epsilon^2 \sqrt{ \sum_{k=1}^{m} \bar b_k^{-4}} + 2\sqrt{2}\sqrt{x_{\beta/2}} \sqrt{ \epsilon^2 \max_{1 \leq l\leq m} (\bar b_l^{-2}) \sum_{k=1}^{m} \left( \frac{b_k}{\bar b_k} \theta_k - \theta_k^0 \right)^2}
$$
which, in turn, holds as soon as
$$
(1-2\sqrt{\delta}\sqrt{x_{\beta/2}}) \sum_{k=1}^{m} \left( \frac{b_k}{\bar b_k} \theta_k - \theta_k^0 \right)^2 \geq t_{\alpha,\bar b} - \epsilon^2 \sum_{k=1}^{m} \bar b_k^{-2} + 2\sqrt{x_{\beta/2}} (1+1/\sqrt{\delta})) \epsilon^2 \sqrt{  \sum_{k=1}^{m} \bar b_k^{-4}}.
$$
By selecting $0<\delta< 1/(4x_{\beta/2})$, and taking into account the definition of $t_{\alpha,\bar b}$, the above inequality holds as soon
$$
\sum_{k=1}^{m} \left( \frac{b_k}{\bar b_k} \theta_k - \theta_k^0 \right)^2 \geq C^{(2)}_{\alpha,\beta} \epsilon^2 \sqrt{\sum_{k=1}^{m} \bar b_k^{-4}},
$$
where
$$
C^{(2)}_{\alpha,\beta} = (c_{\alpha}+2\sqrt{x_{\beta/2}} (1+1/\sqrt{\delta}))/(1-2\sqrt{\delta}\, \sqrt{x_{\beta/2}}).
$$
At this step, we have proved the following assertion
\begin{equation}
\| \bar b^{-1} (b\theta - \bar b\theta^0)\|^2 \geq C^{(2)}_{\alpha,\beta} \rho_{\epsilon,\bar b}^2 \quad \Rightarrow \quad  \mathbb{P}_{\theta,b} \left( \Delta_{\alpha,\bar b}= 0  \right) \leq \frac{\beta}{2}.
\label{eq:imp1}
\end{equation}
Hence, for any $b\in \mathcal{B}$,
$$\mathbb{P}_{\theta,b} (\Delta^{(2)}_{\alpha} = 0) \leq \beta, $$
provided that $\theta\in \ell^2(N_m)$ satisfies
\begin{equation}
 \sum_{k=1}^{m} \bar b_k^{-2}\left( b_k \theta_k - \bar b_k \theta_k^0 \right)^2 \geq C^{(2)}_{\alpha,\beta} \epsilon^2 \sqrt{\sum_{k=1}^{m} \bar b_k^{-4}} \quad \forall\; b, \bar b\in \mathcal{B}.
 \label{eq:cond2}
\end{equation}
Hence, for any $0<\delta< 1/(4x_{\beta/2})$,
$$ \max_{b\in \mathcal{B}} \sup_{\theta\in  \Theta_{1,b}(C_{\alpha,\beta}^{(1)}) \cap \Theta_{2,b}\left(C^{(2)}_{\alpha,\beta}\right)} \mathbb{P}_{\theta,b} (\Delta_\alpha^{(2)}=0)\leq \beta.$$
This completes the proof of the proposition.
\begin{flushright}
$\Box$
\end{flushright}

\subsection{Proof of Proposition \ref{prop:up4}}
\label{s:general_B}

We first investigate the Type I error probability. Remark that that for any $b\in \mathcal{B}$, 
\begin{eqnarray}
 \mathbb{P}_{\theta^0,b} (\Delta_\alpha =1) 
& = &    \mathbb{P}_{\theta^0,b} \left(  \min \left( \Delta_{\alpha,\mathcal{B}_H} , \min_{\tilde b\in \mathcal{B}_{H^c}} \Delta_{\alpha,\tilde b}\right)=1\right) \nonumber \\
& = &    \mathbb{P}_{\theta^0,b} \left( \lbrace  \Delta_{\alpha,\mathcal{B}_H} =1 \rbrace \bigcap_{\tilde b\in \mathcal{B}_{H^c}} \lbrace \Delta_{\alpha,\tilde b}=1\rbrace \right) \nonumber \\
& \leq &   \min \left[ \mathbb{P}_{\theta^0,b} (  \Delta_{\alpha,\mathcal{B}_H} =1 ) , \min_{\tilde b\in \mathcal{B}_{H^c}}  \mathbb{P}_{\theta^0,b}(\Delta_{\alpha,\tilde b} =1)\right].
\label{eq:interB1}
\end{eqnarray}
Remark that a direct consequence of Lemma \ref{lem:1} for $b=\bar b$ is that 
\begin{equation}
\max_{b\in \mathcal{B}}\mathbb{P}_{\theta^0,b}(\Delta_{\gamma, b} =1)\leq \gamma, \quad \forall \gamma \in ]0,1[.
\label{eq:generalB1}
\end{equation}
At this step, we consider two different cases. Assume first that $b \in \mathcal{B}_{H^c}$. Then, starting from (\ref{eq:interB1}),
$$\mathbb{P}_{\theta^0,b} (\Delta_\alpha =1) \leq \min_{\tilde b\in \mathcal{B}_{H^c}}  \mathbb{P}_{\theta^0,b}(\Delta_{\alpha,\tilde b} =1) \leq   \mathbb{P}_{\theta^0,b}(\Delta_{\alpha, b} =1)\leq \alpha,$$
according to (\ref{eq:generalB1}) with $\gamma = \alpha$. Now, assume that $b\in \mathcal{B}_H$. Then, starting again from  (\ref{eq:interB1}), 
\begin{eqnarray*} 
\mathbb{P}_{\theta^0,b} (\Delta_\alpha =1)
&  \leq &  \mathbb{P}_{\theta^0,b} (  \Delta_{\alpha,\mathcal{B}_H} =1 ) \\
& =  & \mathbb{P}_{\theta^0,b} \left(\max_{\tilde b\in \mathcal{B}_H} \Delta_{\frac{\alpha}{|\mathcal{B}_H|}, b}=1\right) \\
& \leq & \sum_{\tilde b\in \mathcal{B}_H}  \mathbb{P}_{\theta^0,b} (\Delta_{\frac{\alpha}{|\mathcal{B}_H|}, \tilde b}=1) \\
& \leq & \alpha,
\end{eqnarray*}
using Lemma \ref{lem:1} since (\ref{eq:cond_reg_gen}) ensures that condition (\ref{eq:cond_reg}) holds for any $b,\tilde b \in \mathcal{B}_H$. The above bounds entail that  
$$ \max_{b\in \mathcal{B}} \mathbb{P}_{\theta^0,b} (\Delta_\alpha =1) \leq \alpha,$$
which ensures that the Type I error probability is controlled by $\alpha$. \\

\noindent
Now, we turn our attention to the control of Type II error probability. For any $\theta\in \ell^2(\mathbb{N}_m)$ and $b\in \mathcal{B}$, we have 
\begin{eqnarray}
\mathbb{P}_{\theta,b} (\Delta_\alpha =0)
& = & \mathbb{P}_{\theta,b} \left(   \min \left( \Delta_{\alpha,\mathcal{B}_H} , \min_{\tilde b\in \mathcal{B}_{H^c}} \Delta_{\alpha,\tilde b}\right) = 0 \right) \nonumber \\
& \leq &   \mathbb{P}_{\theta,b} (   \Delta_{\alpha,\mathcal{B}_H} =0) +  \mathbb{P}_{\theta,b} \left( \min_{\tilde b\in \mathcal{B}_{H^c}} \Delta_{\alpha,\tilde b} = 0 \right) \nonumber \\
& \leq &   \mathbb{P}_{\theta,b} \left(   \max_{\tilde b\in \mathcal{B}_H} \Delta_{\frac{\alpha}{|\mathcal{B}_H|}, \tilde b}=0\right) +  \mathbb{P}_{\theta,b} \left( \min_{\tilde b\in \mathcal{B}_{H^c}} \Delta_{\alpha,\tilde b} = 0 \right) \nonumber \\
& \leq & \min_{\tilde b\in \mathcal{B}_H} \mathbb{P}_{\theta,b} \left(   \Delta_{\frac{\alpha}{|\mathcal{B}_H|}, \tilde b}=0\right) + \sum_{\tilde b \in \mathcal{B}_{H^c}} \mathbb{P}_{\theta,b} \left(  \Delta_{\alpha,\tilde b} = 0 \right). \nonumber
\end{eqnarray}


First, using (\ref{eq:imp1}), 
$$ \min_{\tilde b\in \mathcal{B}_H} \mathbb{P}_{\theta,b} \left(   \Delta_{\frac{\alpha}{|\mathcal{B}_H|}, \tilde b}=0\right) \leq \frac{\beta}{2} \quad \mathrm{provided} \quad
\exists  \, \bar b \in \mathcal{B}_H \quad \mbox{such that} \quad \| \bar b^{-1} (b\theta - \bar b\theta^0)\|^2 \geq C^{(2)}_{\frac{\alpha}{|\mathcal{B}_H|},\beta} \rho_{\epsilon,\bar b}^2, $$
and, hence,
$$\min_{\tilde b\in \mathcal{B}_H} \mathbb{P}_{\theta,b} \left(   \Delta_{\frac{\alpha}{|\mathcal{B}_H|}, \tilde b}=0\right) \leq \frac{\beta}{2} \quad \mathrm{provided} \quad
\theta \in  \bigcup_{\bar b \in \mathcal{B}_H}  \Theta_{3,b}^{\bar b}(C_{\frac{\alpha}{|\mathcal{B}_H|},\beta}).$$
Then, still using (\ref{eq:imp1}), 
$$ \sum_{\tilde b \in \mathcal{B}_{H^c}} \mathbb{P}_{\theta,b} \left(  \Delta_{\alpha,\tilde b} = 0 \right) \leq \frac{\beta}{2} \quad \mathrm{provided} \quad \| \bar b^{-1} (b\theta - \bar b\theta^0)\|^2 \geq C^{(2)}_{\alpha,\frac{\beta}{|\mathcal{B}_H^c|}} \rho_{\epsilon,\bar b}^2 \quad \forall \bar b \in \mathcal{B}_H^c,$$
and, hence,
$$ \sum_{\tilde b \in \mathcal{B}_{H^c}} \mathbb{P}_{\theta,b} \left(  \Delta_{\alpha,\tilde b} = 0 \right) \leq \frac{\beta}{2} \quad \mathrm{provided} \quad \theta \in \bigcap_{\bar b \in \mathcal{B}_H^c}  \Theta_{3,b}^{\bar b}(C_{\alpha,\frac{\beta}{|\mathcal{B}_H^c|}} ).$$
These results holds uniformly in $\theta$, which leads to 
$$  \sup_{\theta\in \Theta_{3,b}} \mathbb{P}_{\theta,b} (\Delta_\alpha=0)\leq \beta.$$
This concludes the proof of the proposition. 

\begin{flushright}
$\Box$
\end{flushright}


\bigskip

\bibliography{Survey}
\bibliographystyle{plain}

\end{document}